\documentclass[12pt,reqno]{amsart}
\usepackage{amsfonts, amsthm, amsmath}
\allowdisplaybreaks[4]
\usepackage{rotating}
\usepackage{tikz}
\usepackage{graphics}
\usepackage{amssymb}
\usepackage{amscd}
\usepackage[latin2]{inputenc}
\usepackage{t1enc}
\usepackage[mathscr]{eucal}
\usepackage{indentfirst}
\usepackage{graphicx}
\usepackage{graphics}
\usepackage{pict2e}
\usepackage{mathrsfs}
\usepackage{enumerate}
\usepackage{CJK}
\usepackage{numprint} 
\usepackage{hyperref} 
\usepackage{cite}
\usepackage{color}
\usepackage{epic}
\usepackage{multirow}
\usepackage{makecell}
\numberwithin{equation}{section}
\theoremstyle{plain}
\newtheorem{theorem}{Theorem}[section]
\newtheorem{lemma}[theorem]{Lemma}
\newtheorem{corollary}[theorem]{Corollary}

\theoremstyle{definition}

\newtheorem{remark}[theorem]{Remark}
\newtheorem{?}[theorem]{Problem}

\begin{document}
\title[Infinitely many negative terms
in two weighted partition sequences]
{Infinitely many negative terms
in two weighted partition sequences}

\author{Meng-Yao Li}
\address[M.-Y. Li]{School of Mathematical Sciences,
	Qufu Normal University, Qufu 273165, P.R. China}
\email{mengyao\_li2026@163.com}

\author{Su-Ping Cui}
\address[S.-P. Cui]{School of Mathematical Sciences, Qufu Normal University, Qufu 273165, P.R. China}
\email{jiayoucui@163.com}

\date{\today}

\begin{abstract}
Recently, Andrews and El Bachraoui introduced partition-counting functions
$a_2(n)$ and $b_2(n)$, together with their signed analogues
$a'_2(n)$ and $b'_2(n)$, and proved that
$a'_2(n), b'_2(n)>0$ for all $n\geq 3$.
For $m\geq 2$, the latter functions extend naturally to
$a'_m(n)$ and $b'_m(n)$, defined by
\begin{align*}
\sum_{n\geq 0} a'_m(n)q^n
&=
-\frac{q}{1-q}
\frac{(q^m;q^m)_\infty}{(q;q^m)_\infty}
+\frac{q}{(1-q)^2},\\
\sum_{n\geq 0} b'_m(n)q^n
&=
-\frac{q(1-q^{m-1})}{1-q}
\frac{(q^{2m};q^m)_\infty}{(q;q^m)_\infty}
+\frac{q}{1-q}
\left(
\frac{1}{1-q}
-\frac{q^{m-1}}{1-q^{m+1}}
\right).
\end{align*}
They asked whether there exists an integer $m>2$
for which $a'_m(n)$ is eventually positive, and posed the
same question for $b'_m(n)$.
In this paper, we resolve these two questions.
\end{abstract}



\maketitle

\section{Introduction}
In this paper, we assume that $|q|<1$ and adopt the standard $q$-series notation
\[
(a;q)_\infty:=\prod_{j=0}^{\infty}(1-aq^j).
\]
The smallest and largest parts of a partition provide natural
ways to refine its enumeration. Partitions subject to restrictions on their smallest and largest parts have been studied extensively. A notable example of a weighted count is Andrews' smallest parts function $\operatorname{spt}(n)$. In \cite{AndrewsSmallestParts}, Andrews derived an identity relating $\operatorname{spt}(n)$,
the partition function $p(n)$, and the second rank moment, and then obtained congruences modulo
$5$, $7$, and $13$. Later, Garvan \cite{Garvan}
proved congruences modulo powers of these primes.
Garvan and Jennings-Shaffer~\cite{GarvanJenningsShaffer} used Bailey-pair techniques to develop a family of two-variable \(\mathrm{spt}\)-crank-type functions. At $z=1$, these functions reduce to the generating functions
of the corresponding smallest-parts functions.
The two-variable functions admit representations in terms
of infinite products or Hecke--Rogers-type double sums.
Restrictions involving both extreme parts lead to a different
class of generating functions. Andrews, Beck, and Robbins
\cite{AndrewsBeckRobbins} considered partitions for which the difference
between the largest and smallest parts is a prescribed
integer $t$. They obtained an explicit rational generating
function when $t>1$ and showed that the corresponding
counting function is a quasipolynomial.
Chern \cite{Chern2017b} studied an overpartition analogue
with bounded differences and related its generating function
to over $q$-binomial coefficients.
Later, Chern and Yee \cite{ChernYee} obtained a further generalization
by giving analytic and combinatorial proofs.
In \cite{Chern2017a}, Chern established a $q$-series identity which
yields a unified derivation of generating functions for
several classes of partitions with bounded differences
between their largest and smallest parts.

When signs are introduced into partition weights, the nonnegativity of the resulting weighted counts becomes a natural question. For weights taking values \(1\) and \(-1\), nonnegativity of the coefficients of the generating function can yield inequalities between partition counts. Andrews and El Bachraoui \cite{AndrewsElBachraouiWeighted} studied several weighted partition generating functions represented by double series. By deriving closed forms for these series and proving that their coefficients are nonnegative, they obtained inequalities between restricted partition counts.

These results motivate the study of how coefficient signs
depend on the parameters in a weighted generating function.
Recently, Andrews and El Bachraoui \cite{AEB} introduced
two partition-counting functions, $a_2(n)$ and $b_2(n)$.
The function $a_2(n)$ counts partitions of $n$ with exactly
$\ell-1$ ones and smallest even part $2k$, such that all
nonunit parts lie in $[2k,2\ell-2]$ and the odd parts greater
than $1$ are distinct.
The function $b_2(n)$ is defined in the same way, with
$[2k,2\ell-2]$ replaced by $[2k,2\ell-1]$.
By assigning weight $(-1)^j$, where
$j$ is the number of odd parts greater than $1$,
they defined the corresponding weighted counts
$a'_2(n)$ and $b'_2(n)$ and proved that
\[
a'_2(n)>0
\qquad\text{and}\qquad
b'_2(n)>0
\qquad(n\ge3).
\]

They further considered the families $a'_m(n)$ and $b'_m(n)$,
defined for integers $m\ge2$ by
\begin{align}
A_m(q)&:=\sum_{n\ge0} a'_m(n)q^n
=
-\frac{q}{1-q}
\frac{(q^m;q^m)_\infty}{(q;q^m)_\infty}
+\frac{q}{(1-q)^2},
\label{eq:A}\\
B_m(q)&:=\sum_{n\ge0} b'_m(n)q^n
=
-\frac{q(1-q^{m-1})}{1-q}
\frac{(q^{2m};q^m)_\infty}{(q;q^m)_\infty}
\nonumber\\
&\qquad
+\frac{q}{1-q}
\left(
\frac{1}{1-q}
-\frac{q^{m-1}}{1-q^{m+1}}
\right).
\label{eq:B}
\end{align}
For $m=3$ and $m=4$, they observed that neither sequence
is positive for all $n\ge3$.
This does not, however, rule out eventual positivity.
They therefore asked whether there exists an integer $m>2$
for which $a'_m(n)$ is eventually positive, and posed the
same question for $b'_m(n)$
\cite[Open Problems 1 and 2]{AEB}.

Our main result answers both questions in the negative.

\begin{theorem}\label{thm:main}
For every integer $m\ge3$, both sequences
$\{a'_m(n)\}_{n\ge0}$ and $\{b'_m(n)\}_{n\ge0}$
contain infinitely many negative terms.
In particular, neither sequence is eventually nonnegative.
\end{theorem}
Throughout, $m\ge3$ is an integer. We use the auxiliary function
\begin{equation}\label{eq:P}
P_m(q):=\frac{(q^m;q^m)_\infty}{(q;q^m)_\infty},
\qquad |q|<1,
\end{equation}
and the notation
\[
a(\theta):=\sum_{k\ge1}\frac{1-\cos(k\theta)}{k^2},
\qquad
b(\theta):=\sum_{k\ge1}\frac{\sin(k\theta)}{k^2},
\qquad 0<\theta<\pi.
\]
In Section~2, we establish the preliminary lemmas.
In Section~3, we prove Theorem~\ref{thm:main}.

\section{Preliminaries}

Let $\mathbb{D}:=\{q\in\mathbb{C}:|q|<1\}$.
We first record an elementary consequence of eventual
nonnegativity.

\begin{lemma}\label{lem:bound}
Let $F(q)=\sum_{n\ge0}c_nq^n$ be analytic in $\mathbb{D}$,
with real coefficients.
If $c_n\ge0$ for all sufficiently large $n$, then there
is a constant $K\ge0$ such that
\begin{equation}\label{eq:circle}
|F(q)|\le F(r)+K
\qquad (|q|=r<1).
\end{equation}
If, in addition,
\[
F(r)\le\frac{r}{(1-r)^2}
\qquad (0<r<1),
\]
then
\begin{equation}\label{eq:bound}
|F(q)|\le\frac{|q|}{(1-|q|)^2}+K
\qquad(q\in\mathbb{D}).
\end{equation}
\end{lemma}

\begin{proof}
Absolute convergence of the Taylor series in $\mathbb{D}$ gives
$$
|F(q)|\le\sum_{n\ge0}|c_n|r^n=\sum_{c_n<0}|c_n|r^n+\sum_{c_n\geq0}|c_n|r^n
=F(r)+2\sum_{c_n<0}|c_n|r^n.
$$

The last sum has finitely many terms and is bounded by $K=2\sum_{c_n<0}|c_n|$. This proves \eqref{eq:circle}.
Then the stated consequence follows immediately, including at $q=0$ by continuity.
\end{proof}

Notice that the products defining $P_m$ converge locally uniformly in $\mathbb{D}$, and neither the numerator nor the denominator vanishes there. In particular, $P_m(q)$, $A_m(q)$ and $B_m(q)$ are analytic in $\mathbb{D}$ and have real Taylor coefficients.
\begin{lemma}\label{lemma2}
For $0<r<1$, we have $P_m(r)>0$.
\end{lemma}
\begin{proof}
We now prove that
\begin{align}
0\leq -\ln(1-x)\leq 2x,
\qquad 0\leq x\leq \frac12. \label{lem2-1}
\end{align}
The left-hand inequality follows from \(0<1-x\leq 1\). For the
right-hand inequality, we use
\[
-\ln(1-x)=\int_0^x\frac{dt}{1-t}.
\]
If \(0\leq x\leq \frac12\), then \(1-t\geq \frac12\) for every
\(t\in[0,x]\). Consequently,
\[
-\ln(1-x)
=
\int_0^x\frac{dt}{1-t}
\leq
\int_0^x 2\,dt
=
2x.
\]

We now consider the numerator. Since \(r^{mj}\to0\) as \(j\to\infty\),
there exists an integer \(J_1\geq1\) such that
\[
r^{mj}\leq \frac12
\qquad (j\geq J_1).
\]
Thus, only finitely many initial terms of the sequence
\(\{r^{mj}\}_{j\geq1}\) need to be treated separately. The corresponding
factors \(1-r^{mj}\) are strictly positive. By \eqref{lem2-1}, for \(j\geq J_1\), we have
\[
0\leq -\ln(1-r^{mj})\leq 2r^{mj}.
\]
Therefore,
\[
\sum_{j=J_1}^{\infty}
\left|\ln(1-r^{mj})\right|
=
\sum_{j=J_1}^{\infty}
\bigl[-\ln(1-r^{mj})\bigr]
\leq
2\sum_{j=J_1}^{\infty}r^{mj}.
\]
The series on the right is geometric, and
\[
\sum_{j=J_1}^{\infty}r^{mj}
=
\frac{r^{mJ_1}}{1-r^m}
<\infty.
\]
It follows that
\[
\sum_{j=1}^{\infty}\ln(1-r^{mj})
\]
converges to a finite real number. Hence
\begin{align}
\prod_{j=1}^{\infty}(1-r^{mj})
=
e^{\left(\sum_{j=1}^{\infty}\ln(1-r^{mj})\right)}>0. \label{lem2-2}
\end{align}
Similarly, we obtain \begin{align}
\prod_{j=0}^{\infty}(1-r^{1+mj})>0. \label{lem2-3}
\end{align}
By definition, $$
P_m(r)
=\frac{\displaystyle\prod_{j=1}^{\infty}(1-r^{mj})}
       {\displaystyle\prod_{j=0}^{\infty}(1-r^{1+mj})}.
$$
Then $P_m(r)>0$ follows from \eqref{lem2-2}, \eqref{lem2-3}, and the above identity.
This completes the proof.
\end{proof}
\begin{corollary}For $0<r<1$, we have
\begin{align}
		A_m(r)&\le\frac{r}{(1-r)^2},\label{A-r}\\
		B_m(r)&\le\frac{r}{(1-r)^2}.\label{B-r}
	\end{align}
\end{corollary}
\begin{proof} Using the fact that
$$
(q^{2m};q^m)_\infty=\frac{(q^m;q^m)_\infty}{1-q^m},
$$
we can rewrite \eqref{eq:A} and \eqref{eq:B} as
\begin{equation}\label{eq:Arewrite}
A_m(q)=\frac{q}{(1-q)^2}
-\frac{q}{1-q}P_m(q)
\end{equation}
and
\begin{equation}\label{eq:Brewrite}
B_m(q)=\frac{q}{(1-q)^2}
-\frac{q^m}{(1-q)(1-q^{m+1})}
-\frac{q(1-q^{m-1})}{(1-q)(1-q^m)}P_m(q).
\end{equation}
Thus, from Lemma \ref{lemma2}, \eqref{eq:Arewrite}, \eqref{eq:Brewrite}, it can be seen that for $m\ge3$ and $0<r<1$,
\begin{equation}\label{eq:real}
	\begin{aligned}
		A_m(r)&=\sum_{n\ge0}a'_m(n)r^n
		=\frac{r}{(1-r)^2}-\frac{r}{1-r}P_m(r)\le\frac{r}{(1-r)^2},\\
		B_m(r)&=\sum_{n\ge0}b'_m(n)r^n
		=\frac{r}{(1-r)^2}-\frac{r^{m}}{(1-r)(1-r^{m+1})}-\frac{r(1-r^{m-1})}{(1-r)(1-r^m)}P_m(r)\nonumber\\
&\quad\le\frac{r}{(1-r)^2}.
	\end{aligned}
\end{equation}
Here we complete the proof.
\end{proof}
We also need the following lemma.
\begin{lemma}\label{lem:trig}
For $0<\theta<\pi$,
the series defining $a(\theta)$
and $b(\theta)$ converge absolutely, and $$a(\theta)>0, \qquad b(\theta)>0.$$
\end{lemma}
\begin{proof}
For $0<\theta<\pi$, we have
\[
\sum_{k\ge1}\left|\frac{1-\cos(k\theta)}{k^2}\right|
=\sum_{k\ge1}\frac{2\sin^2(k\theta/2)}{k^2}
\le 2\sum_{k\ge1}\frac{1}{k^2}<\infty.
\]
Thus the series defining $a(\theta)$ converges absolutely.
Moreover, every term is nonnegative, and the term with
$k=1$ is strictly positive. Hence $a(\theta)>0$.

We next prove that the series defining $b(\theta)$
converges absolutely and that $b(\theta)>0$.

First, the substitution $x=e^{-t}$ gives
\begin{align}
\int_0^1 x^{k-1}(-\ln x)\,dx
=\int_0^\infty te^{-kt}\,dt
=\frac{1}{k^2}. \label{2.10}
\end{align}
Hence $b(\theta)$ can be rewritten as follows.
$$
b(\theta)
=\sum_{k\ge1}\sin(k\theta)
  \int_0^1 x^{k-1}(-\ln x)\,dx.
$$
Since
\begin{align*}
&\quad\sum_{k\ge1}\int_0^1
\left|\sin(k\theta)x^{k-1}(-\ln x)\right|\,dx\\
&=\sum_{k\ge1} \left|\sin(k\theta)\right|
\int_0^1\left|x^{k-1}(-\ln x)\right|\,dx\\
&=\sum_{k\ge1} \left|\sin(k\theta)\right|
\int_0^1x^{k-1}(-\ln x)\,dx\\
&=\sum_{k\ge1}\frac{|\sin(k\theta)|}{k^2}
\leq\sum_{k\ge1}\frac1{k^2}<\infty,
\end{align*}
where we obtain the third step by utilizing \eqref{2.10}. Thus the series defining $b(\theta)$ converges absolutely.

We may interchange summation and integration to obtain
\begin{align}
b(\theta)
=\int_0^1(-\ln x)
  \sum_{k\ge1}\sin(k\theta)x^{k-1}\,dx. \label{3}
\end{align}

For $0<x<1$, the geometric series gives
\begin{align}
\sum_{k\ge1}\sin(k\theta)x^{k-1}
&=\operatorname{Im}\left(
  e^{i\theta}\sum_{j\ge0}(xe^{i\theta})^j
  \right)\nonumber\\
&=\operatorname{Im}\left(
  \frac{e^{i\theta}}{1-xe^{i\theta}}
  \right). \label{1}
\end{align}
Multiplying the numerator and denominator by
$1-xe^{-i\theta}$, we find that
\begin{align}
\frac{e^{i\theta}}{1-xe^{i\theta}}
=\frac{e^{i\theta}-x}{1-2x\cos\theta+x^2}
=\frac{\cos\theta-x+i\sin\theta}
       {1-2x\cos\theta+x^2}. \label{2}
\end{align}
Substituting \eqref{2} into \eqref{1} gives
$$
\sum_{k\ge1}\sin(k\theta)x^{k-1}
=\frac{\sin\theta}{1-2x\cos\theta+x^2}.
$$
Combining \eqref{3} and the above identity, we can obtain that
$$
b(\theta)
=\sin\theta\int_0^1
  \frac{-\ln x}{1-2x\cos\theta+x^2}\,dx.
$$
Finally, $0<\theta<\pi$ implies $\sin\theta>0$, while
for $0<x<1$ we have $-\ln x>0$ and
$$
1-2x\cos\theta+x^2
=(x-\cos\theta)^2+\sin^2\theta>0.
$$
Thus the integrand is strictly positive on $(0,1)$,
and it follows that $b(\theta)>0$.
Here we complete the proof.
\end{proof}
\begin{lemma}\label{prop:growth}
Fix an integer $m\ge3$, and set
\begin{align*}
&\theta=\frac{2\pi}{m},\quad \zeta=e^{i\theta},\quad a=a(\theta),\quad b=b(\theta),\quad
c=\frac{2a}{b},\quad \tau=1+ic.
\end{align*}
For $\varepsilon>0$, let $q_\varepsilon=\zeta e^{-\tau\varepsilon}$. Then $|q_\varepsilon|=e^{-\varepsilon}$ and
\begin{equation}\label{eq:limit}
\lim_{\varepsilon\rightarrow0^+}\varepsilon\ln|P_m(q_\varepsilon)|
=\kappa_m:=\frac{a}{m(1+c^2)}>0.
\end{equation}
In particular, for all sufficiently small $\varepsilon>0$,
\begin{equation}\label{eq:exponential}
|P_m(q_\varepsilon)|\ge e^{\left(\frac{\kappa_m}{2\varepsilon}\right)}.
\end{equation}
\end{lemma}
\begin{proof}
Define the analytic function
\begin{equation}\label{eq:logseries}
L_m(q)=\sum_{k\ge1}\frac{q^k-q^{mk}}{k(1-q^{mk})},\qquad q\in\mathbb{D}.
\end{equation}
For $|q|\le\rho<1$,
\[\sum_{k\ge1}
\left|
\frac{q^k-q^{mk}}{k(1-q^{mk})}
\right|
\le
\sum_{k\ge1}\frac{\rho^k+\rho^{mk}}{k(1-\rho^m)}.
\]
The series on the right converges, so the series defining
$L_m$ converges locally uniformly in $\mathbb D$.
Hence $L_m$ is analytic in $\mathbb D$.
We now show that \begin{equation}\label{eq:realpart}
\ln|P_m(q)|=\operatorname{Re}L_m(q).
\end{equation}

Recall that, for $|z|<1$, the principal logarithm satisfies
\[
-\operatorname{Ln}(1-z)
=
\sum_{k\geq1}\frac{z^k}{k}.
\]
Taking real parts gives
\[
\ln|1-z|
=
-\operatorname{Re}
\sum_{k\geq1}\frac{z^k}{k}.
\]
Applying this identity to the factors in the product representation of
\(P_m(q)\), we obtain
\begin{align*}
\ln|P_m(q)|&=\ln\left|
\frac{\displaystyle\prod_{j=1}^{\infty}(1-q^{mj})}
     {\displaystyle\prod_{j=0}^{\infty}(1-q^{1+mj})}
\right|\\
&=
\sum_{j\geq1}\ln|1-q^{mj}|
-
\sum_{j\geq0}\ln|1-q^{1+mj}|\\
&=
-\operatorname{Re}
\sum_{j\geq1}\sum_{k\geq1}
\frac{q^{mjk}}{k}
+
\operatorname{Re}
\sum_{j\geq0}\sum_{k\geq1}
\frac{q^{(1+mj)k}}{k}\nonumber\\
&=
\operatorname{Re}
\sum_{k\geq1}\frac{1}{k}
\left(
-\sum_{j\geq1}q^{mjk}
+
\sum_{j\geq0}q^{(1+mj)k}
\right)\nonumber\\
&=
\operatorname{Re}
\sum_{k\geq1}
\frac{1}{k}
\left(
-\frac{q^{mk}}{1-q^{mk}}
+
\frac{q^k}{1-q^{mk}}
\right)\\
&=
\operatorname{Re}
\sum_{k\geq1}
\frac{q^k-q^{mk}}{k(1-q^{mk})}\\
&=
\operatorname{Re}L_m(q).
\end{align*}

Since $\zeta^m=1$ and $q_\varepsilon=\zeta e^{-\tau\varepsilon}$, we have
\begin{equation}\label{eq:scaledlog}
\varepsilon L_m(q_\varepsilon)
=\sum_{k\ge1}\frac{\varepsilon(\zeta^ke^{-k\tau\varepsilon}-e^{-mk\tau\varepsilon})}
{k(1-e^{-mk\tau\varepsilon})}.
\end{equation}
We now prove that $\varepsilon L_m(q_\varepsilon)$
converges as $\varepsilon\to0^+$. Using $\operatorname{Re}\tau=\operatorname{Re}(1+ic)=1$ and $m\ge1$, we obtain
\begin{align}
\left|\frac{\varepsilon(\zeta^ke^{-k\tau\varepsilon}-e^{-mk\tau\varepsilon})}
{k(1-e^{-mk\tau\varepsilon})}\right|
&\le\frac{2\varepsilon e^{-k\varepsilon}}{k(1-e^{-mk\varepsilon})}\nonumber\\
&\le\frac{2\varepsilon}{k(e^{k\varepsilon}-1)}\le\frac2{k^2}.
\label{eq:majorant}
\end{align}
The majorant is summable and independent of $\varepsilon$. Dominated convergence therefore yields
\begin{align*}
\lim_{\varepsilon\rightarrow0^+}\varepsilon L_m(q_\varepsilon)
&=\frac1{m\tau}\sum_{k\ge1}\frac{\zeta^k-1}{k^2}\nonumber\\
&=\frac1{m\tau}\sum_{k\ge1}\frac{e^{ik\theta}-1}{k^2}\nonumber\\
&=\frac1{m\tau}\sum_{k\ge1}\frac{\cos(k\theta)+i\sin(k\theta)-1}{k^2}\nonumber\\
&=\frac{-a+ib}{m(1+ic)}.
\end{align*}
Taking real parts and utilizing \eqref{eq:realpart} together
with the identity $bc=2a$, we obtain
$$
\lim_{\varepsilon\rightarrow0^+}\varepsilon\ln|P_m(q_\varepsilon)|
=\frac{-a+bc}{m(1+c^2)}=\frac{a}{m(1+c^2)}>0.
$$
This proves \eqref{eq:limit}.
Since the limit in \eqref{eq:limit} is positive, we find that
for all sufficiently small $\varepsilon>0$,
 $$\varepsilon\ln|P_m(q_\varepsilon)|\geq\frac{a}{2m(1+c^2)}.$$
Namely, $$|P_m(q_\varepsilon)|\geq e^{\frac{a}{2m\varepsilon(1+c^2)}}.$$
Here we complete the proof of \eqref{eq:exponential}.
\end{proof}
\begin{remark}
The direction of approach is essential. Along the radial
path $q_\varepsilon=\zeta e^{-\varepsilon}$, the same
calculation gives
\[
\lim_{\varepsilon\to0^+}
\varepsilon\log|P_m(q_\varepsilon)|=-\frac{a}{m}<0.
\]
The choice $\tau=1+ic$ with $c=2a/b$ makes this limit
positive while preserving
$|q_\varepsilon|=e^{-\varepsilon}<1$.
\end{remark}

\section{Proof of the main theorem}

\begin{proof}[Proof of Theorem~\ref{thm:main}]
Fix an integer $m\geq3$ and use Lemma \ref{prop:growth}. Thus
$$
\zeta=e^{2\pi i/m},\qquad
\tau=1+ic,\qquad
q_\varepsilon=\zeta e^{-\tau\varepsilon}.
$$
Since $\operatorname{Re}\tau=1$, we have
$$
|q_\varepsilon|=e^{-\varepsilon}<1,
\qquad
q_\varepsilon\longrightarrow\zeta
\quad(\varepsilon\rightarrow0^+).
$$
Put $\alpha=\kappa_m/2>0$. By \eqref{eq:exponential},
there exists $\varepsilon_0>0$ such that
\begin{align}\label{pmm}
|P_m(q_\varepsilon)|\geq e^{\alpha/\varepsilon}
\qquad(0<\varepsilon<\varepsilon_0).
\end{align}
All constants below may depend on $m$, which is fixed,
but are independent of $\varepsilon$.

We consider the two sequences separately.

\medskip
\noindent
(i) The sequence $\{a'_m(n)\}_{n\ge0}$.
Suppose, to the contrary, that $a'_m(n)$ has only finitely
many negative terms. Then $a'_m(n)\geq0$ for all sufficiently
large $n$. By Lemma~\ref{lem:bound} and \eqref{A-r},
there exists a constant $K_A\geq0$ such that
$$
|A_m(q_\varepsilon)|
\leq A_m(e^{-\varepsilon})+K_A
\leq \frac{e^{-\varepsilon}}
           {(1-e^{-\varepsilon})^2}+K_A.
$$
For $0<\varepsilon\leq1$, we have
$$
1-e^{-\varepsilon}
=\int_0^\varepsilon e^{-t}\,dt
\geq e^{-1}\varepsilon.
$$
Since $e^{-\varepsilon}\leq1$, it follows that
\begin{equation}\label{eq:Aupper}
|A_m(q_\varepsilon)|
\leq e^2\varepsilon^{-2}+K_A
\leq C_A\varepsilon^{-2},
\qquad 0<\varepsilon\leq1,
\end{equation}
where $C_A=e^2+K_A$.

On the other hand, \eqref{eq:A} gives
\begin{align}
A_m(q_\varepsilon)
=\frac{q_\varepsilon}{(1-q_\varepsilon)^2}
-\frac{q_\varepsilon}{1-q_\varepsilon}
 P_m(q_\varepsilon). \label{Am-1}
\end{align}
Since $\zeta\neq1$, we have
\begin{align*}
\lim_{\varepsilon\rightarrow0^+}\frac{q_\varepsilon}{1-q_\varepsilon}
=\lim_{\varepsilon\rightarrow0^+}\frac{\zeta e^{-\tau \varepsilon}}{1-\zeta e^{-\tau \varepsilon}}
=\frac{\zeta}{1-\zeta}\neq0
\end{align*}
and
\begin{align*}
\lim_{\varepsilon\rightarrow0^+}\frac{q_\varepsilon}{(1-q_\varepsilon)^2}
=\lim_{\varepsilon\rightarrow0^+}\frac{\zeta e^{-\tau \varepsilon}}{(1-\zeta e^{-\tau \varepsilon})^2}
=
\frac{\zeta}{(1-\zeta)^2}.
\end{align*}
Consequently, for all sufficiently small $\varepsilon>0$,
$$
\left|\frac{q_\varepsilon}{1-q_\varepsilon}\right|
\geq d_A:=
\frac12\left|\frac{\zeta}{1-\zeta}\right|>0,
$$
and
$$
\left|\frac{q_\varepsilon}{(1-q_\varepsilon)^2}\right|
\leq M_A,
\qquad
M_A:=1+\left|\frac{\zeta}{(1-\zeta)^2}\right|.
$$
Combining \eqref{pmm}, \eqref{Am-1}, and the two inequalities
above, we obtain, by the reverse triangle inequality,
\begin{align*}
|A_m(q_\varepsilon)|
&\geq
\left|\frac{q_\varepsilon}{1-q_\varepsilon}\right|
|P_m(q_\varepsilon)|
-
\left|\frac{q_\varepsilon}{(1-q_\varepsilon)^2}\right|\\
&\geq d_Ae^{\alpha/\varepsilon}-M_A.
\end{align*}
From \eqref{eq:Aupper} and the above inequality, it can be seen that
$$
d_Ae^{\alpha/\varepsilon}-M_A
\leq C_A\varepsilon^{-2}.
$$
Multiplying by $\varepsilon^2$ gives
$$
d_A\varepsilon^2e^{\alpha/\varepsilon}
\leq C_A+M_A\varepsilon^2.
$$
The right-hand side tends to $C_A$, whereas the left-hand
side tends to infinity. Indeed, with $t=1/\varepsilon$,
$$
\varepsilon^2e^{\alpha/\varepsilon}
=\frac{e^{\alpha t}}{t^2}
\geq\frac{(\alpha t)^3}{6t^2}
=\frac{\alpha^3t}{6}
\longrightarrow\infty.
$$
This contradiction shows that $a'_m(n)$ has infinitely
many negative terms.

\medskip
\noindent
(ii) The sequence $\{b'_m(n)\}_{n\ge0}$.
Suppose now that $b'_m(n)$ has only finitely many negative
terms. Applying Lemma~\ref{lem:bound} and \eqref{B-r}
as above, we obtain a constant $K_B\geq0$ such that
$$
|B_m(q_\varepsilon)|
\leq B_m(e^{-\varepsilon})+K_B
\leq\frac{e^{-\varepsilon}}
          {(1-e^{-\varepsilon})^2}+K_B.
$$
Thus, with $C_B=e^2+K_B$,
\begin{align}
|B_m(q_\varepsilon)|
\leq C_B\varepsilon^{-2}
\qquad(0<\varepsilon\leq1). \label{Bm-1}
\end{align}

To obtain a lower bound, write \eqref{eq:Brewrite} as
\begin{align}
B_m(q)=R_m(q)-H_m(q)P_m(q), \label{Bm}
\end{align}
where
$$
R_m(q)
=\frac{q}{(1-q)^2}
-\frac{q^m}{(1-q)(1-q^{m+1})},
\qquad
H_m(q)
=\frac{q(1-q^{m-1})}{(1-q)(1-q^m)}.
$$
Since $\zeta^m=1$ and $\zeta^{m+1}=\zeta\neq1$,
\begin{align*}
\lim_{\varepsilon\rightarrow0^+}R_m(q_\varepsilon)
&=\lim_{\varepsilon\rightarrow0^+}\left(\frac{\zeta e^{-\tau \varepsilon }}{(1-\zeta e^{-\tau \varepsilon })^2}
-\frac{e^{-m\tau \varepsilon }}{(1-\zeta e^{-\tau \varepsilon })(1-\zeta e^{-(m+1)\tau \varepsilon })}\right)\nonumber\\
&=
\frac{\zeta}{(1-\zeta)^2}
-\frac{1}{(1-\zeta)^2}
=-\frac{1}{1-\zeta}.
\end{align*}
In particular, for all sufficiently small $\varepsilon>0$,
\begin{align}
|R_m(q_\varepsilon)|\leq M_B,
\qquad
M_B:=1+\frac{1}{|1-\zeta|}. \label{Rm}
\end{align}

Next, since $q_\varepsilon^m=e^{-m\tau\varepsilon}$,
the Taylor expansion of the exponential gives
$$
1-q_\varepsilon^m
=1-e^{-m\tau\varepsilon}
=m\tau\varepsilon+O(\varepsilon^2).
$$
As $\tau\neq0$, it follows that
\begin{align*}
\lim_{\varepsilon\rightarrow0^+}\frac{\varepsilon}{1-q_\varepsilon^m}
=\frac{1}{m\tau}.
\end{align*}
Utilizing the above identity, we subsequently arrive at
\begin{align}\label{eq:Hlimit}
\lim_{\varepsilon\rightarrow0^+}
\varepsilon H_m(q_\varepsilon)
&=
\lim_{\varepsilon\rightarrow0^+}
\left(
\frac{q_\varepsilon(1-q_\varepsilon^{m-1})}
     {1-q_\varepsilon}
\frac{\varepsilon}{1-q_\varepsilon^m}
\right)\nonumber\\
&=
\lim_{\varepsilon\to0^+}
\left(
\frac{\zeta e^{-\tau\varepsilon}
      \left(1-\zeta^{m-1}e^{-(m-1)\tau\varepsilon}\right)}
     {1-\zeta e^{-\tau\varepsilon}}
\cdot
\frac{\varepsilon}{1-q_\varepsilon^m}
\right)\nonumber\\
&=\frac{\zeta(1-\zeta^{-1})}{(1-\zeta)m\tau}\nonumber\\
&=-\frac{1}{m\tau}\neq0.
\end{align}

Taking absolute values in \eqref{eq:Hlimit}, we obtain
$$
\lim_{\varepsilon\rightarrow0^+}\varepsilon|H_m(q_\varepsilon)|
=\frac{1}{m|\tau|}>0.
$$
Hence, for all sufficiently small $\varepsilon>0$,
$$
|H_m(q_\varepsilon)|
\geq\frac{d_B}{\varepsilon},
\qquad
d_B:=\frac{1}{2m|\tau|}>0.
$$

Then in view of \eqref{pmm}, \eqref{Bm}, \eqref{Rm}, the above inequality, and the reverse triangle inequality, we deduce that
\begin{align*}
|B_m(q_\varepsilon)|
&\geq
|H_m(q_\varepsilon)|\,|P_m(q_\varepsilon)|
-|R_m(q_\varepsilon)|\\
&\geq
\frac{d_B}{\varepsilon}e^{\alpha/\varepsilon}-M_B.
\end{align*}
Combining this with \eqref{Bm-1} gives
$$
\frac{d_B}{\varepsilon}e^{\alpha/\varepsilon}-M_B
\leq C_B\varepsilon^{-2},
$$
or equivalently,
$$
d_B\varepsilon e^{\alpha/\varepsilon}
\leq C_B+M_B\varepsilon^2.
$$
This is again impossible. Indeed, setting $t=1/\varepsilon$,
we have
$$
\varepsilon e^{\alpha/\varepsilon}
=\frac{e^{\alpha t}}{t}
\geq\frac{(\alpha t)^2}{2t}
=\frac{\alpha^2t}{2}
\longrightarrow\infty,
$$
whereas $C_B+M_B\varepsilon^2\to C_B$.
Thus $b'_m(n)$ also has infinitely many negative terms.

Since $m\geq3$ was arbitrary, the theorem follows.
\end{proof}

\subsection*{Acknowledgements}
The second author was supported by the Taishan Scholar Project of Shandong Province (No.tsqn202507173), the Natural Science Foundation of Shandong Province of China (No. ZR2025MS91).

\subsection*{Declaration of AI usage.}
 During the preparation of this work, the authors used AI to explore proof
strategies and to assist with language editing and proofreading. The authors take full responsibility for the final content of the publication.

\end{document}